\documentclass[a4paper,12pt]{article}

\usepackage{amsmath,amsfonts,amssymb,amsthm}
\theoremstyle{plain}
\newtheorem{thm}{Theorem}
\newtheorem*{thm*}{Theorem}
\newtheorem*{Mthm}{Main Theorem}
\newtheorem*{TECHthm}{Technical Theorem}
\newtheorem*{CLT}{Local CLT}
\newtheorem{lem}{Lemma}
\newtheorem*{lem*}{Lemma}
\newtheorem*{cor}{Corollary}
\newtheorem*{prop}{Proposition}

\theoremstyle{definition}

\newtheorem{Def}{Definition}

\newcommand{\reftit}{\textit}    
\newcommand{\refis}{\textbf}     

\begin{document}

\title{Brownian Bridge Asymptotics for the Subcritical Bernoulli
 Bond Percolation.}

\author{Yevgeniy Kovchegov\\
 \small{Email: yevgeniy@math.stanford.edu}\\
 \small{Fax: 1-650-725-4066}}

\maketitle

\begin{abstract}
 For the $d$-dimensional model of a subcritical bond percolation ($p<p_c$) and a
 point $\mathbf{\vec{a}}$ in $\mathbb{Z}^d$,
 we prove that
a cluster conditioned on connecting points $(0,...,0)$ and $n
\mathbf{\vec{a}}$ if scaled by $\frac{1}{n \| \mathbf{\vec{a}}
\|}$ along $\mathbf{\vec{a}}$ and by $\frac{1}{\sqrt{n}}$ in the
orthogonal direction converges asymptotically to Time $\times$
($d-1$)-dimensional Brownian Bridge.
\end{abstract}

\section{Introduction.}

\subsection{Percolation and Brownian Bridge.}

 We begin by briefly stating the notion of a bond percolation based
 on the material rigorously presented in \cite{grimmett} and \cite{kesten}, and
 the notion of the Brownian Bridge as well as the word description
 of the result connecting the two that we have obtained and
 made the primary objective of this paper.

 \textbf{Percolation:} For each edge of the $d$-dimensional square lattice
 $\mathbb{Z}^d$ in turn, we declare the edge $open$ with probability
 $p$ and $closed$ with probability $1-p$, independently of all
 other edges. If we delete the closed edges, we are left with a
 random subgraph of $\mathbb{Z}^d$. A connected component of the
 subgraph is called a ``cluster", and the number of edges in a
 cluster is the ``size" of the cluster. The probability $\theta (p)$ that
 the point $(0,0)$ belongs to a cluster of an infinite size is zero if
 $p=0$, and one if $p=1$. However, there exists a critical
 probability $0< p_c <1$ such that $\theta (p) = 0$ if $p < p_c$
 and $\theta (p) >0$ if $p > p_c$. In the first case, we say that
 we are dealing with a $subcritical$ percolation model, and in
 the second case, we say that we are dealing with a
 $supercritical$ percolation model.

 \textbf{Brownian Bridge:} defined as a sample-continuous Gaussian
 process $B^0$ on $[0,1]$ with mean $0$ and $\bold{E} B^0_s B^0_t
 = s(1-t)$ for $0 \leq s \leq t \leq 1$. So, $B^0_0 = B^0_1 = 0$
 a.s. Also, if $B$ is a Brownian motion, then the process
 $B_t - tB_1$ ($0 \leq t \leq 1$) is a Brownian Bridge. For more
 details see \cite{bill}, \cite{dudley} and \cite{durrett}. The
 process $B^{0, \mathbf{\vec{a}}}_t \equiv B^0_t +t \mathbf{\vec{a}}$
 is a Brownian Bridge connecting points zero and
 $\vec{\mathbf{a}}$.

 \textbf{History of the problem:}  Below, we consider the $d$-dimensional model
 of a subcritical bond percolation ($p<p_c$) and a point $\mathbf{\vec{a}}$
 in $\mathbb{Z}^d$, conditioned on the event of zero being connected to
 $n \mathbf{\vec{a}}$. We first show that a specifically
 chosen path connecting points zero and $n \mathbf{\vec{a}}$ and going
 through some appropriately defined points on the cluster (regeneration
 points), if scaled $\frac{1}{n \| \mathbf{\vec{a}}\|}$ times along
 $\mathbf{\vec{a}}$ and $\frac{1}{\sqrt{n}}$ times in the
 direction orthogonal to $\mathbf{\vec{a}}$, converges to Time $\times$
 ($d-1$)-dimensional Brownian Bridge as $n \rightarrow +\infty$,
 where the scaled interval connecting points zero and $n \mathbf{\vec{a}}$
 serves as a $[0,1]$ time interval. In other words, we prove that a
 scaled ``skeleton" going through the regeneration points of the
 cluster converges to Time $\times$ ($d-1$)-dimensional Brownian
 Bridge. In a subsequent step, we show that if scaled, then the hitting
 area of the orthogonal hyper-planes shrinks, implying that for
 $n$ large enough, all the points of the scaled cluster are within
 an $\varepsilon$-neighborhood of the points in the ``skeleton".
 One of the major tools used in this research was the renewal
 technique developed in \cite{acc}, \cite{ccc}, \cite{ioffe}, \cite{cc}
 and \cite{ioffe1}  as part of the derivation of the Ornstein-Zernike estimate
 for the subcritical bond percolation model and ``similar" processes.
 A major result related to the study is  that for
 $\mathbf{\vec{a}}=(1,0,...,0)$,
 the hitting distribution of the cluster in the
 intermediate planes, $x_1 =tn\mathbf{\vec{a}}$, $0<t<1$ obeys a
 multidimensional local limit theorem (see \cite{ccc}). Dealing with all other
 $\mathbf{\vec{a}} \not= (k,0,...,0)$ became possible only after the
 corresponding technique further mastering the regeneration structures and
 equi-decay profiles was developed in \cite{ioffe} and \cite{ioffe1}.
 This technique played a central role in obtaining the research results.

\subsection{Asymptotic Convergence.}\label{intro:2}

 Here we state a version of a local CLT and a technical result that we
 later prove.\\

 \textbf{Local Limit Theorem:} In this paper we are going to use the version of
 the local CLT borrowed from \cite{durrett}:
 Let $X_1, X_2,... \in \mathbb{R}$ be i.i.d. with $\bold{E} X_i =0$,
 $\bold{E} X_i^2 = \sigma ^2 \in (0, \infty)$, and having
 a common lattice distribution with span $h$. If
 $S_n = X_1 + ... + X_n$ and $P[X_i \in b + h \mathbb{Z}] = 1$
 then $P[S_n \in nb + h \mathbb{Z}] = 1$. We put
 $$p_n(x) = P[S_n / \sqrt{n} =x] \mbox{  for } x \in \Lambda_n
   = \{ (nb + hz) / \sqrt{n} \mbox{ : } z \in \mathbb{Z} \} $$
 and
 $$n(x)= (2 \pi \sigma^2)^{-1/2} \exp(-x^2/2\sigma^2) \mbox{  for }
   x \in (-\infty, \infty)$$

 \begin{CLT}
 Under the above hypotheses,
 $\sup_{x \in \Lambda_n} |{\frac{\sqrt{n}}{h} p_n(x) - n(x)}|
  \rightarrow 0$ as $n \rightarrow \infty$.
 \end{CLT}

 \textbf{Technical Result Concerning Convergence to the Brownian Bridge}
 (to be used in Chapter 2.4, and is proved in Chapter 3.2):
 The following technical result is going to be proved in the section \ref{bb} of
 this paper, however since we are to use it in the section \ref{percol},
 we will state the result below as part of the introduction.
 Let $X_1, X_2,...$ be i.i.d. random variables on $\mathbb{Z}^d$
  with the span of the lattice distribution equal to one
(see \cite{durrett}, section 2.5), and let there be a
 $\bar{\lambda} > 0$ such that the moment-generating function
 $$\bold{E}(e^{\theta \cdot X_1}) <\infty$$
for all $\theta \in B_{\bar{\lambda}}$.
\\

 Now, for a given vector $\mathbf{\vec{a}} \in \mathbb{Z}^d$,
 let $X_1+...+X_i =[t_i, Y_i]_f \in \mathbb{Z}^d$ when
 written in the new orthonormal basis such that
 $\mathbf{\vec{a}} = [\| \mathbf{\vec{a}} \|, 0]_f$ (in the new basis
 $[\cdot , \cdot ]_f \in \mathbb{R} \times \mathbb{R}^{d-1}$).
 Also let $P[\mathbf{\vec{a}} \cdot X_i]>0]=1$.
 We define the process $[t, Y_{n,k}^*(t)]_f$ to be
 the interpolation of $0$ and
 $[\frac{1}{n \| \mathbf{\vec{a}} \|}t_i, \frac{1}{\sqrt{n}}Y_i]_f^{i=0,1,...,k}$,
 in Section 2.2 we will show that

\begin{TECHthm}
The process
$$\{ Y^*_{n,k}\mbox{ for some } k \mbox{ such that } [t_k, Y_k]_f =n \mathbf{\vec{a}} \}$$
conditioned on the existence of such $k$ converges  weakly to the
 Brownian Bridge (of variance  that depends only on the law of $X_1$).
\end{TECHthm}

\section{The Main Result in Subcritical Percolation.}
\label{percol}

In this section we work only with subcritical percolation
probabilities $p<p_c$.

\subsection{Preliminaries.}

Here we briefly go over the definitions that one can find in
Section 4 of \cite{ioffe}.
\\
We start with the inverse correlation length $\xi_p(\vec{x})$:

$$\xi_p(\vec{x}) \equiv -\lim_{n \rightarrow \infty}
   \frac{1}{n} P_p (0 \leftrightarrow [n\vec{x}]),$$
where the limit is always defined due to the FKG property of the
Bernoulli bond percolation (see \cite{grimmett}).
 Now, $\xi_p(\vec{x})$ is the support function of the compact convex set
$$ \bold{K}^p \equiv \bigcap_{\vec{n} \in \mathbb{S}^{d-1}}
  \lbrace \vec{r} \in \mathbb{R}^d \mbox{ : }
  \vec{r} \cdot \vec{n} \leq \xi_p(\vec{n}) \rbrace ,$$
with non-empty interior int\{$\bold{K}^p$\} containing point zero.
\\
Let $\mathbf{\vec{r}} \in \partial \bold{K}^p$, and let $\vec{e}$
be a basis vector such that $\vec{e} \cdot \mathbf{\vec{r}}$ is
maximal. For $\vec{x},\vec{y} \in \mathbb{Z}^d$ define
$$S^r_{\vec{x},\vec{y}} \equiv \lbrace \vec{z} \in \mathbb{R}^d \vert
  \mathbf{\vec{r}} \cdot \vec{x} \leq \mathbf{\vec{r}} \cdot \vec{z}
  \leq \mathbf{\vec{r}} \cdot \vec{y} \rbrace .$$
Note that $S^r_{\vec{x},\vec{y}} = \emptyset$ if
 $ \mathbf{\vec{r}} \cdot \vec{y} < \mathbf{\vec{r}} \cdot \vec{x} $.
\\
Let $\bold{C}^r_{\vec{x},\vec{y}}$ denote the corresponding common
open cluster of $x$ and $y$ when we run the percolation process on
$S^r_{\vec{x},\vec{y}}$.

\begin{Def}
For $\vec{x},\vec{y} \in \mathbb{Z}^d$ lets define
$h_r$-connectivity
 $\lbrace \vec{x} \leftarrow^{h_r} \rightarrow \vec{y} \rbrace$
 of $\vec{x}$ and $\vec{y}$ to be the event that\\
1. $\vec{x}$ and $\vec{y}$ are connected in the restriction of the
percolation configuration to the slab $S^r_{\vec{x},\vec{y}}$.
\\
2. If $\vec{x} \not= \vec{y}$, then
  $\bold{C}^r_{\vec{x},\vec{y}} \bigcap S^r_{\vec{x},\vec{x}+\vec{e}}
   = \lbrace \vec{x},\vec{x}+\vec{e} \rbrace$
and
  $\bold{C}^r_{\vec{x},\vec{y}} \bigcap S^r_{\vec{y}-\vec{e},\vec{y}}
   = \lbrace \vec{y}-\vec{e},\vec{y} \rbrace$.
\\
3. If $\vec{x}=\vec{y}$ and all the edges adjoined to $\vec{x}$
and perpendicular to $\vec{e}$ are closed.
\end{Def}

 Set
 $$h_r(\vec{x}) \equiv P_p \lbrack 0 \leftarrow^{h_r} \rightarrow \vec{x} \rbrack .$$
Notice that $h_r(0) = (1-p)^{2(d-1)}$.

\begin{Def}
For $\vec{x},\vec{y} \in \mathbb{Z}^d$ lets define
$f_r$-connectivity
 $\lbrace \vec{x} \leftarrow^{f_r} \rightarrow \vec{y} \rbrace$
 of $\vec{x}$ and $\vec{y}$ to be the event that\\
 1. $\vec{x} \not= \vec{y}$\\
 2. $\vec{x} \leftarrow^{h_r} \rightarrow \vec{y}$ .\\
 3. For no $\vec{z} \in \mathbb{Z}^d \setminus \lbrace \vec{x},\vec{y} \rbrace $  both
$ \lbrace \vec{x} \leftarrow^{h_r} \rightarrow \vec{z} \rbrace
\mbox{ and  }
   \lbrace \vec{z} \leftarrow^{h_r} \rightarrow \vec{y} \rbrace$
take place.
\end{Def}

Set
 $$f_r(\vec{x}) \equiv P_p \lbrack 0 \leftarrow^{f_r} \rightarrow \vec{x} \rbrack .$$
Notice that $f_r(0)=0$.

\begin{Def}
Suppose $0 \leftarrow^{h_r} \rightarrow \vec{x}$, we say that
$\vec{z} \in \mathbb{Z}^d$ is \textbf{a regeneration point}
 of $\bold{C}_{0,\vec{x}}^r$ if
\\
 1. $\mathbf{\vec{r}} \cdot \vec{e} \leq \mathbf{\vec{r}} \cdot \vec{z}
     \leq \mathbf{\vec{r}} \cdot (\vec{y}-\vec{e}) $ \\
 2. $S_{\vec{z}-\vec{e},\vec{z}+\vec{e}}^r \bigcap \bold{C}_{0,\vec{x}}^r$
    contains exactly three points: $\vec{z}-\vec{e}$, $\vec{z}$ and
    $\vec{z}+\vec{e}$, where $\vec{e}$ is defined as before.\\
Let also $\vec{x}$ itself be a regeneration point.
\end{Def}

The following Ornstein-Zernike equality is due to be used soon:
\begin{thm*}
 $\exists$ $A(\cdot, \cdot)$ on $(0,p_c) \times \mathbf{S}^{d-1}$ s. t.
 \begin{eqnarray} \label{oz}
 P_p[0 \leftrightarrow \vec{x}] = \frac{A(p, n(\vec{x}))}
 {\|\vec{x}\|^{\frac{d-1}{2}} } e^{-\xi_p(\vec{x})} (1+o(1))
 \end{eqnarray}
 uniformly in $\vec{x} \in \mathbb{Z}^d$,
 where $n(\vec{x}) \equiv \frac{\vec{x}}{\|\vec{x}\|}$.
\end{thm*}
We refer to \cite{ioffe} for the proof of the theorem.

\subsection{Measure $Q_{r_0}^r(x)$.}

It had been proved in section 4 of \cite{ioffe} that for a given
$\mathbf{\vec{r}}_0 \in \partial \bold{K}^p$ there exists
$\bar{\lambda} >0$ such that
 $$F_{r_0}(\mathbf{\vec{r}})=
  \frac{1}{(1-p)^{2(d-1)}} \sum_{x \in \mathbb{Z}^d} f_{\mathbf{\vec{r}}_0}(x)
  e^{\mathbf{\vec{r}} \cdot \vec{x}} =1 \mbox{ whenever }
 \mathbf{\vec{r}} \in B_{\bar{\lambda}}(\mathbf{\vec{r}}_0)
 \bigcap \partial \bold{K}^p$$
and therefore
$$Q_{r_0}^r(\vec{x}) \equiv \frac{1}{(1-p)^{2(d-1)}}f_{r_0}(\vec{x})
  e^{\mathbf{\vec{r}} \cdot \vec{x}}
 \mbox{  is a measure on  } \mathbb{Z}^d .$$
 Also, it was shown that
$$\mu=\mu_{r_0}(\mathbf{\vec{r}}) \equiv \bold{E}_{r_0}^rX =
 \sum_{\vec{x} \in \mathbb{Z}^d}\vec{x}Q_{r_0}^r(\vec{x})
 = \nabla_r logF_{r_0}(\mathbf{\vec{r}}) \not= 0$$
and
 $$F_{r_0}(\mathbf{\vec{r}}) < \infty \mbox{  for all } \mathbf{\vec{r}} \mbox{ in }
   B_{\bar{\lambda}}(\mathbf{\vec{r}}_0).$$

The later implies
$$F_{r_0}(\mathbf{\vec{r}}) =\sum_{\vec{x} \in \mathbb{Z}^d}f_{r_0}(\vec{x})
  e^{\mathbf{\vec{r}} \cdot \vec{x}}
  = \sum_{\vec{x} \in \mathbb{Z}^d}Q_{r_0}^{r_0}(\vec{x}) e^{\theta \cdot \vec{x}}
 < \infty$$
 for $\theta = \mathbf{\vec{r}}-\mathbf{\vec{r}}_0 \in B_{\bar{\lambda}}(0)$,\\
i.e. the moment generating function
 $\bold{E}_{r_0}^{r_0} (e^{\theta \cdot X_1})$
 of the law $Q_{r_0}^{r_0}$ is finite for all
 $\theta \in B_{\bar{\lambda}}(0)$.\\

Now, there is a renewal relation (see section 1 and section 4 of
\cite{ioffe}),
 $$ h_{r_0}(\vec{x})=\frac{1}{(1-p)^{2(d-1)}}
    \sum_{\vec{z} \in \mathbb{Z}^d} f_{r_0}(\vec{z})h_{r_0}(\vec{x}-\vec{z})
  \mbox{  with  } h_{r_0}(0)=(1-p)^{2(d-1)}$$
and therefore
$$h_{r_0}([N \mu])=(1-p)^{2(d-1)}e^{-r \cdot [N \mu]}
     \sum_{k} \bigotimes_1^k Q_{r_0}^r(X_1+...+X_k=[N \mu])
  \mbox{  for  } N>0 ,$$

where $X_1,X_2,...$ is a sequence of i.i.d. random variables
distributed according to $Q_{r_0}^r$, as $h_{r_0}$-connection
is a chain of $f_{r_0}$-connections with junctions
at the regeneration points of $\bold{C}_{0,x}^{r_0}$.
\\

\subsection{Important Observation.}

The probability that $0 \leftarrow^{h_{r_0}} \rightarrow x$ with
 exactly $k$ regeneration points $x_1, x_1+x_2, ... , \sum_{i=1}^k x_i =x$
\begin{eqnarray} \label{imp}
  P_X & \equiv &
  P[0 \leftarrow^{h_{r_0}} \rightarrow x \mbox{ ;  regeneration points:  }
   x_1, x_1+x_2, ... , \sum_{i=1}^k x_i =x] \nonumber \\
  & = & \frac{1}{(1-p)^{2(d-1)(k-1)}}
  P[0 \leftarrow^{f_{r_0}} \rightarrow x_1]
  P[x_1  \leftarrow^{f_{r_0}} \rightarrow x_1+x_2]...
  P[\sum_{i=1}^{k-1} x_i
     \leftarrow^{f_{r_0}} \rightarrow \sum_{i=1}^k x_i =x] \nonumber \\
 & = & \frac{1}{(1-p)^{2(d-1)(k-1)}} f_{r_0}(x_1)f_{r_0}(x_2)...f_{r_0}(x_k).
\end{eqnarray}

\subsection{The Result.}

In this section we fix $\bold{\vec{a}} \in \mathbb{Z}^d$, and let
 $\bold{r}=\bold{r}_0= \bold{\vec{a}} \mathbb{R}^+\bigcap \partial \bold{K}^p$.
Then we recall that
 $$\bold{E}_{r_0}^r (e^{\theta \cdot X_1}) < \infty$$
 for all $\theta \in B_{\bar{\lambda}}(0)$.
 We also denote $h(x) \equiv h_{r_0}(x)$ and $f(x) \equiv f_{r_0}(x)$.
 \\
First, we introduce a new basis $\{ \vec{f_1},\vec{f_2},...,
\vec{f_d} \}$, where
 $\vec{f_1} = \frac{\bold{\vec{a}}}{\| \bold{\vec{a}} \|}$. We use
 $[\cdot,\cdot]_f \in \mathbb{R} \times \mathbb{R}^{d-1}$ to denote
 the coordinates of a vector with respect to the new basis.
 Obviously $\mathbf{\vec{a}}=[\|\mathbf{\vec{a}} \|, 0]_f$.
 We want to prove that the process corresponding to the last $d-1$
coordinates in the new basis of the scaled
 ($\frac{1}{n \| \bold{\vec{a}} \|}$ times along $\bold{\vec{a}}$
 and $\frac{1}{\sqrt{n}}$ times in the orthogonal d-1
 dimensions) interpolation of regeneration points
  of $\bold{C}_{0, n \bold{\vec{a}}}^{r_0}$  conditioned on
 ${0 \leftarrow^{h}\rightarrow n \bold{\vec{a}}}$
 converges weakly to the Brownian Bridge $B^o(t)$ (with variance that depends
 only on measure $Q_{r_0}^r$) where $t$ represents the scaled first coordinate
 in the new basis.\\

Let $X_1, X_2,...$ be i.i.d. random variables distributed
according to $Q_{r_0}^r$ law.
 We interpolate $0,X_1,(X_1+X_2),...,(X_1+...+X_k)$ and scale by
$\frac{1}{n\|\mathbf{\vec{a}} \|} \times \frac{1}{\sqrt{n}}$ along
$<\mathbf{\vec{a}}> \times <\mathbf{\vec{a}}>^{\bot}$ to get the
process $[t, Y^*_{n,k}(t)]_f$. The technical theorem (see Chapters
(\ref{intro:2}) and (\ref{bb:gen})) implies the following

\begin{thm}
The process
$$\{ Y^*_{n,k}\mbox{ for some } k \mbox{ such that }
    X_1+...+X_k= n \bold{\vec{a}} \}$$
conditioned on  the existence of such $k$
 converges weakly to the
Brownian Bridge (with variance that depends
only on measure $Q_{r_0}^r$).
\end{thm}

Now, let for $y_1,...,y_k \in \mathbb{Z}^d$ with  positive
increasing first coordinates
 $\gamma (y_1,...,y_k)$  be the last $(d-1)$ coordinates
  in the new basis
 of the scaled ($\frac{1}{n \|\mathbf{\vec{a}}\|} \times \frac{1}{\sqrt{n}}$)
interpolation of points $0,y_1,...,y_k$ (where the first
coordinate is time). Notice that $\gamma (y_1,...,y_k) \in
C_o[0,1]^{d-1}$
 as a function of scaled first coordinate
 whenever $y_k= n \bold{\vec{a}}$.
\\
 By the important observation (\ref{imp}) we've made before,
 for any function $F(\cdot )$ on $C[0,1]^{d-1}$,\\
\\
$\sum_k \sum_{x_1+...+x_k= n \bold{\vec{a}}}
  F(\gamma (x_1, x_1+x_2, ... , \sum_{i=1}^k x_i))$
$$ \times
 P[0 \leftarrow^{h_{r_0}} \rightarrow x
  \mbox{ ;  regeneration points:  }
  x_1, x_1+x_2, ... , \sum_{i=1}^k x_i =x]$$

$$=\sum_k \sum_{x_1+...+x_k= n \bold{\vec{a}}}
 F(\gamma (x_1, x_1+x_2, ... , \sum_{i=1}^k x_i))
 \frac{1}{(1-p)^{2(d-1)(k-1)}} f(x_1)...f(x_k)$$

$$= (1-p)^{2(d-1)} e^{-r \cdot n \bold{\vec{a}}}
 \sum_k \sum_{x_1+...+x_k= n \bold{\vec{a}}}
 F(\gamma (x_1, x_1+x_2, ... , \sum_{i=1}^k x_i))
 Q_{r_0}^r(x_1)...Q_{r_0}^r(x_k).$$
\\
Therefore, for any $A \subset C[0,1]^{d-1}$\\
\\
$P_p[ \gamma (\mbox{regeneration points}) \in A \mbox{  } | \mbox{  }
      0 \leftarrow^{h}\rightarrow  n \bold{\vec{a}} ]$

$$
=\frac{  \sum_k \sum_{x_1+...+x_k= n \bold{\vec{a}}}
         I_A (\gamma (x_1, x_1+x_2, ... , \sum_{i=1}^k x_i))
         \frac{1}{(1-p)^{2(d-1)(k-1)}} f(x_1)...f(x_k) }
        { \sum_k \sum_{x_1+...+x_k= n \bold{\vec{a}}}
          \frac{1}{(1-p)^{2(d-1)(k-1)}} f(x_1)...f(x_k) }$$

$$
=\frac{\sum_k \sum_{x_1+...+x_k= n \bold{\vec{a}}}
         I_A (\gamma (x_1, x_1+x_2, ... , \sum_{i=1}^k x_i))
         Q_{r_0}^r(x_1)...Q_{r_0}^r(x_k)}
        {\sum_k \sum_{x_1+...+x_k= n \bold{\vec{a}}}
         Q_{r_0}^r(x_1)...Q_{r_0}^r(x_k)}$$
$$= P[Y^*_{n,k} \in A \mbox{ for the } k \mbox{ such that }
    X_1+...+X_k= n \bold{\vec{a}} \mbox{  } | \mbox{  }
  \exists k \mbox{ such that }
    X_1+...+X_k= n \bold{\vec{a}}] .$$

Hence, we have proved the following

\begin{cor}
 The process corresponding to the last $d-1$ coordinates (in the new basis
 $\{ \vec{f_1},\vec{f_2},...,\vec{f_d} \}$) of the scaled
 $({\frac{1}{n \|\mathbf{\vec{a}} \|} \times \frac{1}{\sqrt{n}} })$
interpolation of regeneration points
  of $\bold{C}_{0,n \bold{\vec{a}}}^{r_0}$ (where the first coordinate
 is time) conditioned on
${0 \leftarrow^{h}\rightarrow n \bold{\vec{a}}}$
 converges
weakly to the Brownian Bridge (with variance that depends
only on measure $Q_{r_0}^r$).
\end{cor}


\subsection{Shrinking of the Cluster. Main Theorem.}

Here for $\bold{\vec{a}} \in \mathbb{Z}^d$ we let
 $\bold{r}_0= \bold{\vec{a}} \mathbb{R}^+\bigcap \partial \bold{K}^p$ again.
Before we proceed with the proof that the scaled percolation
cluster $\bold{C}_{0,n \bold{\vec{a}}}^{r_0}$ shrinks to the
scaled interpolation skeleton of regeneration points, we need to
prove the following
\begin{prop}
If $\mathbf{\vec{r}} = \nabla \xi_p(\mathbf{\vec{r}}_0)$ then
$Q_{r_0}^r$ is a probability measure.
\end{prop}

\begin{proof}
 First we notice that $\mathbf{\vec{r}}_0 \cdot \mathbf{\vec{r}}
 = \mathbf{\vec{r}}_0 \cdot \nabla \xi_p(\mathbf{\vec{r}}_0)
 = D_{\mathbf{\vec{r}}_0}(\xi_p(\mathbf{\vec{r}}_0)) =
 \xi_p(\mathbf{\vec{r}}_0)$, and thus
$$H_{r_0}(\mathbf{\vec{r}}) \equiv
 \frac{1}{(1-p)^{2(d-1)}} \sum_{\vec{x} \in \mathbb{Z}^d}
  h_{r_0}(x)e^{\mathbf{\vec{r}} \cdot \vec{x}} \geq
 \sum_{\vec{x} \in <\mathbf{\vec{a}}> \cap \mathbb{Z}^d}
  h_{r_0}(x)e^{\mathbf{\vec{r}} \cdot \vec{x}}
 = \sum_{\vec{x} \in <\mathbf{\vec{a}}> \cap \mathbb{Z}^d}
    h_{r_0}(x)e^{\xi_p(\vec{x})}
 = +\infty$$
for $d \leq 3$ by Ornstein-Zernike equation (\ref{oz}). For all
other $d$ we sum over all $\vec{x}$ inside a small enough cone
around
$\mathbf{\vec{a}}$ to get $H_{r_0}(\mathbf{\vec{r}}) = +\infty$.\\

 Now, for all $\vec{n} \in \mathbb{S}^{d-1}$, $\vec{n} \cdot
 \nabla \xi_p(\mathbf{\vec{r}}_0) = D_{\vec{n}} \xi_p(\mathbf{\vec{r}}_0)
 \leq \xi_p(\vec{n})$ by convexity of $\xi_p$, and therefore
 $\mathbf{\vec{r}} = \nabla \xi_p(\mathbf{\vec{r}}_0) \in
  \partial \mathbf{K}^p$. Notice that due to the strict convexity of
  $\xi_p$ and the way $\mathbf{K}^p$ was defined,
  $\mathbf{\vec{r}} = \nabla \xi_p(\mathbf{\vec{r}}_0)$ is the
  only point on $\partial \mathbf{K}^p$ such that
  $\mathbf{\vec{r}}_0 \cdot \mathbf{\vec{r}} = \xi_p(\mathbf{\vec{r}}_0)$.\\

 Now, Ornstein-Zernike equation (\ref{oz}) also implies that the sums
 $H_{r_0}(\tilde{r})$ and $F_{r_0}(\tilde{r})$
 are finite whenever
 $\tilde{r} \in \alpha \bold{K}^p =
  \bigcap_{\vec{n} \in \mathbb{S}^{d-1}}
  \lbrace \vec{r} \in \mathbb{R}^d \mbox{ : }
  \vec{r} \cdot \vec{n} \leq \alpha \xi_p(\vec{n}) \rbrace $
  with $\alpha \in (0,1)$, and
 due to the recurrence relation of $f_{r_0}$ and $h_{r_0}$
 connectivity functions,
 $H_{r_0}(\tilde{r}) = \frac{1}{1 -F_{r_0}(\tilde{r})}$
 (see \cite{ioffe}). Therefore
 $F_{r_0}(\mathbf{\vec{r}}) \equiv
 \frac{1}{(1-p)^{2(d-1)}} \sum_{\vec{x} \in \mathbb{Z}^d}
  f_{r_0}(x)e^{\mathbf{\vec{r}} \cdot \vec{x}} = 1$,
  where the probability measure $Q_{r_0}^r$ has an
  exponentially decaying tail due to the same reasoning as in chapter 4
  of \cite{ioffe} ("mass-gap" property).

\end{proof}

 With the help of the proposition above we shell show that the consequent
 regeneration points are situated relatively close to each other:
\begin{lem*}
$$P_p[\max_{i} |x_i - x_{i-1}|>n^{1/3},
 \mbox{  } x_i \mbox{- reg. points  }
  | \mbox{  } 0 \leftarrow^h \rightarrow n \bold{\vec{a}} ]
 <\frac{1}{n}$$
for $n$ large enough.
\end{lem*}

\begin{proof}
 Let $\mathbf{\vec{r}} \equiv \nabla \xi_p(\mathbf{\vec{r}}_0)
 = \nabla \xi_p(\mathbf{\vec{a}})$.
 Since $\xi_p(x)$ is strictly convex
 (see section 4 in \cite{ioffe}),
$$ \frac{ \xi_p(\bold{\vec{a}}) - \xi_p(\bold{\vec{a}} -\frac{\vec{x}}{n})}
 {(\frac{\|\vec{x}\|}{n})}
 < \frac{\vec{x}}{\|\vec{x}\|} \cdot \nabla \xi_p(\bold{\vec{a}}) $$
for $\vec{x} \in \mathbb{Z}^d$ ($\vec{x} \not= 0$), and therefore
$$\xi_p(n \bold{\vec{a}})- \xi_p(n \bold{\vec{a}}-\vec{x})
 = \|\vec{x}\| \frac{ \xi_p(\bold{\vec{a}}) - \xi_p(\bold{\vec{a}}
-\frac{\vec{x}}{n})}{(\frac{\|\vec{x}\|}{n})}
 < \vec{x} \cdot \nabla \xi_p(\bold{\vec{a}})
 = \mathbf{\vec{r}} \cdot \vec{x}.$$
Thus, since $Q_{r_0}^r(x)$ decays exponentially and therefore
$$\frac{f(x)}{(1-p)^{2(d-1)}}
 e^{\xi_p(n \bold{\vec{a}})- \xi_p(n \bold{\vec{a}}-x)} < Q_{r_0}^r(x)$$
and also decays exponentially. Hence by Ornstein-Zernike result
(\ref{oz}),

$$P_p[  n^{1/3} <
  |x| , \mbox{ } x \mbox{-first reg. point  }
  | 0 \leftarrow^h \rightarrow n \bold{\vec{a}} ]
 =\sum_{ n^{1/3} < |x| }
 \frac{f(x)}{(1-p)^{2(d-1)}} \frac{h(n \bold{\vec{a}}-x)}
      {h(n \bold{\vec{a}})}
 < \frac{1}{n^2} $$
for $n$ large enough.
 So, since the number of the regeneration points is no greater
 than $n$,

 $$P_p[\max_{i} |x_i - x_{i-1}|>n^{1/3},
 \mbox{  } x_i \mbox{- reg. points  }
  | \mbox{  } 0 \leftarrow^h \rightarrow n \bold{\vec{a}} ]
 <\frac{1}{n}$$
for $n$ large enough.

\end{proof}

 Now,  it is really easy to check that there is a
constant $\lambda_f >0$ such that
 $$f(\vec{x}) > e^{-\lambda_f \|\vec{x}\|}$$
for all $\vec{x}$ such that $f(\vec{x}) \not= 0$ (here we only
need to connect points $\vec{e}$ and $\vec{x} - \vec{e}$ with two
non-intersecting open paths surrounded by the closed edges), and
there exists a constant $\lambda_u >0$ such that
 $$P_p[\mbox{ percolation cluster } \bold{C}(0)
  \not\subset [\mathbb{R}; B_{R}^{d-1}(0)]_f ] < e^{-\lambda_u R} $$
for $R$ large enough due to the exponential decay of the radius
distribution for subcritical probabilities (see \cite{grimmett}).
Hence, for a given $\epsilon
>0$
$$P_p[\mbox{ cluster } \bold{C}_{0,\vec{x}}^{r_0}
 \not\subset [\mathbb{R}, B_{\epsilon \sqrt{n}}^{d-1}(0)]_f
 \mbox{  } | \mbox{  } 0 \leftarrow^f \rightarrow x]
  < e^{ \lambda_f \|\vec{x}\| -\lambda_u \epsilon \sqrt{n}}, $$
and therefore, summing over the regeneration points, we get

$$P_p[\mbox{ scaled cluster } \bold{C}_{0, n \bold{\vec{a}} }^{r_0}
 \not\subset \epsilon \mbox{-neighbd. of }
 [0,1] \times \gamma (\mbox{ reg. points })
 \mbox{  } | \mbox{  } 0 \leftarrow^g \rightarrow  n \bold{\vec{a}} ]$$

$$ < \frac{1}{n} +
 n e^{ \lambda_f n^{1/3} -\lambda_u \epsilon \sqrt{n}} $$
 for $n$ large enough.\\

We can now state the main result of this paper:
\begin{Mthm}
 The process corresponding to the last $d-1$ coordinates (in the new basis
 $\{ \vec{f_1},\vec{f_2},...,\vec{f_d} \}$) of the scaled
 $({\frac{1}{n \|\mathbf{\vec{a}} \|} \times \frac{1}{\sqrt{n}} })$
interpolation of regeneration points
  of $\bold{C}_{0,n \bold{\vec{a}}}^{r_0}$ (where the first coordinate
 is time) conditioned on
${0 \leftarrow^{h}\rightarrow n \bold{\vec{a}}}$
 converges
weakly to the Brownian Bridge (with variance that depends only on
measure $Q_{r_0}^r$).
\\
 Also for a given $\epsilon >0$
$$P_p[\mbox{ scaled cluster } \bold{C}_{0, n \bold{\vec{a}} }^{r_0}
 \not\subset \epsilon \mbox{-neighbd. of }
 [0,1] \times \gamma (\mbox{ reg. points })
 \mbox{  } | \mbox{  } 0 \leftarrow^h \rightarrow  n \bold{\vec{a}} ]
 \rightarrow 0$$
as $n \rightarrow \infty$.
\end{Mthm}

\section{Convergence to Brownian Bridge.} \label{bb}

As it was mentioned in the introduction, this chapter is entirely
dedicated to proving the Technical Theorem that we have already
used in the proof of the main result.

\subsection{Simple Case.}

 Let $Z_1, Z_2,...$ be i.i.d. random variables on $\mathbb{Z}$ with
the span of the lattice distribution equal to one (see
\cite{durrett}, section 2.5) and mean $\mu = \bold{E}Z_1 <\infty$,
$\sigma^2=Var(Z_1)<\infty$. Also let point zero be inside of the
closed convex hull of $\{ z \mbox{ : } P[Z_1 = z]>0 \}$.
\\

Consider a one dimensional plane and a walk $X_j$ that starts with
$X_0=0$ and for a given $X_j$, the (j+1)-st step  to be
$X_{j+1}=X_j + Z_{j+1}$. After interpolation we get
$$X(t)=X_{[t]}+(t-[t])(X_{[t]+1}-X_{[t]})$$
for $0\leq{t}<\infty$.\\
And define $\bar{X}(t) = (t,X(t))$ to be a two dimensional walk.\\

Now, if for a given integer $n>0$  we define
$X_n(t)\equiv{\frac{X(nt)}{\sqrt{n}}}$ for $0\leq{t}\leq{1}$,
 then
$X_n(t)$ would belong to $C[0,1]$ and $X_n(0)=0$.

\begin{thm}\label{simpleT}
$X_n(t)$ conditioned on $X_n(1)=0$ converges weakly to the
Brownian Bridge.
\end{thm}

First we need to prove the theorem when $\mu =0$. For this we need
to prove that

\begin{lem}
For $A_0 \subseteq {C[0,1]}$, let
$P_n(A_0)=P[X_n\in{A_0}|X_n(1)=0]$ to be the law of $X_n$
conditioned on $X_n(1)=0$.
Then \\

(a) For $\mu=0$, the finite-dimensional distributions of $P_n$
converge weakly to a Gaussian distributions.\\

(b) There are positive $\{C_n\}_{n=1,2,...}\rightarrow{C}$ ($C
=\sigma^2$ when $\mu=0$) such that $0<C<\infty$ and
$$Cov_{P_n}(X_n(s),X_n(t)) = C_ns(1-t) + O(\frac{1}{n})$$
for all $0\leq{s}\leq{t}\leq{1}$.  More precisely:
$Cov_{P_n}(X_n(s),X_n(t))=C_ns(1-t)$ if $[ns]<[nt]$ and \\
$Cov_{P_n}(X_n(s),X_n(t))=C_ns(1-t) - C_n\frac{\epsilon_1 (1-
\epsilon_2 )}{n}$ if $[ns]=[nt]$, where
 $\epsilon_1 = \frac{ns-[ns]}{n} \in [0,1)$ and
 $\epsilon_2 = \frac{nt-[nt]}{n} \in [0,1)$.
\end{lem}

and we need

\begin{lem}
For $\mu=0$, the probability measures $P_n$ induced on the
subspace of $X_n(t)$ trajectories in $C[0,1]$ are tight.
\end{lem}

\begin{proof}[Proof of Lemma 1:]
(a) Though it is not difficult to show that a finite-dimensional
distribution of $P_n$ converges weakly to a gaussian distribution,
here we only show the convergence for one and two points on the
interval (in case of one point $t \in [0,1]$, we show that the
limit variance has to be equal to $t(1-t) \sigma^2$). Take $t \in
\frac{1}{n}\mathbb{Z} \cap (0,1)$ and let $\alpha
=\frac{k}{\sqrt{n}}$, then by the Local CLT,

 \begin{eqnarray} \label{phi}
P[X(tn)=k] = \frac{1}{\sqrt{n}} \Phi_{\sigma \sqrt{t}}(\alpha) +
o(\frac{1}{\sqrt{n}}), \mbox{  where  } \Phi_{v}(x) \equiv
\frac{1}{v \sqrt{2\pi}}e^{-\frac{x^2}{2 v^2}}
 \end{eqnarray} is
the normal density function, and the error term is uniformly
bounded by a $o(\frac{1}{\sqrt{n}})$ function independent of $k$.
\\
Therefore, substituting (\ref{phi}),
$$P_n[X_n(t)=\alpha]  =  \frac{(\frac{1}{\sqrt{n}} \Phi_{\sigma
\sqrt{t}}(\alpha) + o(\frac{1}{\sqrt{n}}))
 (\frac{1}{\sqrt{n}} \Phi_{\sigma \sqrt{1-t}}(\alpha)+ o(\frac{1}{\sqrt{n}}))}
 {\frac{1}{\sqrt{n}} \Phi_{\sigma}(0) + o(\frac{1}{\sqrt{n}}) }
 = \frac{1}{\sqrt{n}} \Phi_{\sigma \sqrt{t(1-t)}}(\alpha)
 + o(\frac{1}{\sqrt{n}}) .$$
Thus for a set $A$ in $\mathbb{R}$,
$$P_n[X_n(t)\in{A}] = \sum_{k\in\sqrt{n}A}[\frac{1}{\sqrt{n}}
\Phi_{\sigma \sqrt{t(1-t)}}(\alpha)+ o(\frac{1}{\sqrt{n}})] =
N[0, t(1-t)\sigma^2](A) + o(1)$$
 -here the limit variance is equal to $t(1-t)\sigma^2$. Given that
the variance $\sigma^2 <0$, the convergence follows.
\\

The same method works for more than one point, here we do it for
two: Let $\alpha_1 = \frac{k_1}{\sqrt{n}}$ and $\alpha_2 =
\frac{k_2}{\sqrt{n}}$,
 then as before, for $t_1<t_2$ in $\frac{1}{n}\mathbb{Z} \cap (0,1)$,
 writing the conditional probability as a ratio of two probabilities,
 and representing the probabilities according to (\ref{phi}),
 we get
 $$P_n[X_n(s)=\alpha_1 , X_n(t)=\alpha_2] =
 \frac{\sqrt{|\mathcal{A}|}}{2\pi \sigma^2 }
   \exp \left \{-\frac{(\alpha_1, \alpha_2)\mathcal{A} (\alpha_1, \alpha_2)^T}
{2\sigma^2} \right \} +o(\frac{1}{n}) .
$$
\\
\\
where
$$\mathcal{A}={\left( \begin{array}{cc}
             \frac{t_2}{(t_2-t_1)t_1} & -\frac{1}{t_2-t_1} \\
             -\frac{1}{t_2-t_1} & \frac{1-t_1}{(t_2-t_1)(1-t_2)} \end{array} \right)}.$$
\\
Thus for sets $A_1$ and $A_2$ in $\mathbb{R}$,

\begin{eqnarray*}
P_n[X_n(t_1) \in A_1 , X_n(t_2) \in A_2] & = &
 \sum_{k_1 \in \sqrt{n}A_1, k_2 \in \sqrt{n}A_2}
 [\frac{\sqrt{|\mathcal{A}|}}{2\pi \sigma^2 }
  \exp \left \{ -\frac{(\alpha_1, \alpha_2) \mathcal{A}  (\alpha_1, \alpha_2)^T}
    {2\sigma^2} \right \}+o(\frac{1}{n})] \\
\\
& = & N[0, \mathcal{A}^{-1}](A_1 \times A_2) + o(1)
\end{eqnarray*}
 Observe that
 $(\sigma^2 \mathcal{A}^{-1})={\left( \begin{array}{cc}
                  t_1(1-t_1)\sigma^2 & t_1(1-t_2)\sigma^2 \\
                  t_1(1-t_2)\sigma^2 & t_2(1-t_2)\sigma^2 \end{array} \right)}$
is the covariance matrix, and the part (b) of the lemma follows in case $\mu =0$.\\

(b) Though the estimate above produces the needed variance in case
when the mean $\mu =0$ , in general, we need to apply the
following approach: We first consider the case when $s<t$ and both
$s,t \in \frac{1}{n}\mathbb{Z} \cap (0,1)$ where
$$\mathbf{E}[X_n(s) \mbox{  }|\mbox{  } X_n(t) =y]
  =  \mathbf{E}[Z_1 +...+ Z_{sn} |Z_1 +...+ Z_{tn} = y]
  =  \frac{s}{t} y ,$$
and therefore

\begin{eqnarray*}
Cov_{P_n}(X_n(s),X_n(t)) & = &
 \frac{s}{t} \mathbf{E}[X^2_n(t)|X_n(1)=0]
\end{eqnarray*}
as $\{ -X_n(1-t) \mbox{  }|\mbox{  } X_n(1)=0 \}$ and $\{ X_n(t)
\mbox{  }|\mbox{  } X_n(1)=0 \}$ are identically distributed.\\

Now, by symmetry (time reversal),
 $$Cov_{P_n}(X_n(s),X_n(t)) =Cov_{P_n}(X_n(1-t),X_n(1-s))
 =\frac{1-t}{1-s} \mathbf{E}[X^2_n(s)|X_n(1)=0],$$
 and therefore
$$\frac{\mathbf{E}[X^2_n(s)|X_n(1)=0]}
  {\mathbf{E}[X^2_n(t)|X_n(1)=0]}=\frac{s(1-s)}{t(1-t)}.$$
Hence, there exists a constant $C_n$ such that for all $t \in
\frac{1}{n}\mathbb{Z} \cap (0,1)$
$$\frac{\mathbf{E}[X^2_n(t)|X_n(1)=0]}{t(1-t)} \equiv C_n.$$
Thus we have shown that for $s \leq t$ in $\frac{1}{n}\mathbb{Z}
\cap [0,1]$,
$$Cov_{P_n}(X_n(s),X_n(t))  =
\frac{s}{t} \mathbf{E}[X^2_n(t)|X_n(1)=0]
  =  \frac{s}{t}C_n t(1-t)
  =  C_n s(1-t).$$

Now, consider the general case: $s=s_0 + \frac{\epsilon_1}{n} \leq
t = t_0 + \frac{\epsilon_2}{n}$, where $ns_0, nt_0 \in \mathbb{Z}$
and $\epsilon_1, \epsilon_2 \in [0,1)$. Then the covariance


\begin{eqnarray*}
 Cov_{P_n}(X_n(s),X_n(t))
 & = & (1-\epsilon_1)(1-\epsilon_2)Cov_{P_n}(X_n(s_0),X_n(t_0))\\
 & + & (1-\epsilon_1)\epsilon_2Cov_{P_n}(X_n(s_0),{X_n(t_0+\frac{1}{n})}) \\
 & + &\epsilon_1(1-\epsilon_2)Cov_{P_n}(X_n(s_0+\frac{1}{n}),X_n(t_0))\\
 & + &\epsilon_1\epsilon_2Cov_{P_n}(X_n(s_0+\frac{1}{n}),X_n(t_0+\frac{1}{n}))
\end{eqnarray*}
Therefore
\begin{eqnarray*}
 Cov_{P_n}(X_n(s),X_n(t))
 & = & C_ns(1-t) \mbox{ when  } s_0<t_0 \mbox{  (}[ns]<[nt]\mbox{),}
\end{eqnarray*}
and
\begin{eqnarray*}
 Cov_{P_n}(X_n(s),X_n(t))
 & = & C_ns(1-t) - C_n\frac{\epsilon_1 (1- \epsilon_2 )}{n}
  \mbox{ when  } s_0=t_0 \mbox{  (}[ns]=[nt]\mbox{).}
\end{eqnarray*}
Now, plugging in $s=t=\frac{1}{2}$ we get
$$C_n =4\mathbf{E}[X^2_n(\frac{1}{2})|X_n(1)=0] \mbox{  when }n\mbox{  is even,}$$
and
$$C_n =4\mathbf{E}[X^2_n(\frac{1}{2})|X_n(1)=0] (\frac{n}{n-1})
  \mbox{  when }n\mbox{  is odd.}$$
Therefore
$$C_n =4\mathbf{E}[X^2_n(\frac{1}{2})|X_n(1)=0]
      (1+O(\frac{1}{n})) \rightarrow{C} =\sigma^2$$
as $\{ X_n(\frac{1}{2}), P_n \}$ converges in distribution as $n
\rightarrow +\infty$.

\end{proof}

\begin{proof}[Proof of Lemma 2:]

 Before we begin the proof of tightness, we notice that the only
 real obstacle we face is that the process is conditioned on
 $X_n=0$. The tightness for the case without the conditioning has been proved years
 ago as part of the Donsker's Theorem (see Chapter 10 in
 \cite{bill}). With the help of the local CLT we are essentially
 removing the difference between the two cases.

 Given a $\lambda >0$ and let $m=[n \delta]$ for
 a given $0< \delta \leq 1$, then for any $\mu >0$,

\begin{eqnarray*}
 P_{\lambda} & \equiv & P[ \max_{0\leq i \leq m} X_i \geq
 \lambda \sqrt{n} > X_m> -\lambda \sqrt{n}| X_n = 0] \\
 & = & \sum_{a=-[\lambda \sqrt{n}]}^{[\lambda \sqrt{n}]}
  \frac{P[\max_{0 \leq i \leq m}X_i > [\lambda \sqrt{n}]
         \mbox{ ; }  X_m =a
         \mbox{ ; } X_n=0]}
       {P[X_n = 0]}\\
 & = &
 \sum_{a=-[\lambda \sqrt{n}]}^{[\lambda \sqrt{n}]}
  \frac{P[\max_{0 \leq i \leq m}X_i > [\lambda \sqrt{n}]
         \mbox{ ; }  X_m =a]
         P[X_{n-m}=-a]}
       {P[X_n = 0]}\\
 & \leq &
 \max_{-[\lambda \sqrt{n}] \leq a \leq [\lambda \sqrt{n}]}
  (\frac{P[X_{n-m}=-a]}
       {P[X_n = 0]})
 \times
  \sum_{a=-[\lambda \sqrt{n}]}^{[\lambda \sqrt{n}]}
  P[\max_{0 \leq i \leq m}X_i > [\lambda \sqrt{n}]
         \mbox{ ; }  X_m =a]
 \\
 & \leq &
 2P[\max_{0\leq i \leq m}X_i \geq \lambda \sqrt{n} \geq X_m \geq
     -\lambda \sqrt{n} ]
\end{eqnarray*}
\\
for $n$ large enough, where by the local CLT,

$$ \max_{-[\lambda \sqrt{n}] \leq a \leq [\lambda \sqrt{n}]}
  (\frac{P[X_{n-m}=-a]}
       {P[X_n = 0]})
\leq 2
$$
for $n$ large enough as $n-m$ linearly depends on $n$.
\\

 Therefore, the probability
 $$P[\max_{0\leq i \leq m} |X_i| \geq \lambda \sqrt{n} |X_n=0]
   \leq 2P_{\lambda} + P[|X_m| \geq \lambda \sqrt{n} |X_n=0],$$

 where
 $$P_{\lambda}
 \leq 2P[\max_{0\leq i \leq m} X_i \geq \lambda \sqrt{n}].$$

Now, due to the point-wise convergence, we can proceed as in
Chapter 10 of \cite{bill} by bounding the two remaining
probabilities:

$$P[\max_{0\leq i \leq m} |X_i| \geq \lambda \sqrt{n}]
 \leq 2 P[|X_m| \geq \frac{1}{2} \lambda \sqrt{n} ] \rightarrow
 2P[|\sqrt{\delta}N| \geq \frac{\lambda}{2 \sigma}]
\leq \frac{16 \delta^{3/2} \sigma^3}{\lambda^3}
\mathbf{E}[|N|^3]$$ and similarly
$$P[|X_m| \geq \lambda \sqrt{n} | X_n=0] \rightarrow
 P[|\sqrt{\delta(1-\delta)}N| \geq \frac{\lambda}{\sigma}]
\leq \frac{\delta^{3/2} \sigma^3}{\lambda^3} \mathbf{E}[|N|^3].$$

Thus, for all integer $k \in [0, n-m]$,
$$P[\max_{0\leq i \leq m} |X_{k+i} - X_k| \geq \lambda \sqrt{n}
|X_n=0] = P[\max_{0\leq i \leq m} |X_i| \geq \lambda \sqrt{n}
|X_n=0] \leq 70\frac{\delta^{3/2} \sigma^3}{\lambda^3}
\mathbf{E}[|N|^3]$$ for $n$ large enough, (see Chapter 10 in
\cite{bill}). Therefore $\{ P_n \}$ are tight (see Chapter 8 of
\cite{bill}).

\end{proof}

\begin{proof}[Proof of Theorem \ref{simpleT}:]
The lemmas above imply the convergence when the mean $\mu =0$.
Now, for $\mu \not= 0$, there exists a $\rho \in \mathbb{R}$ such
that $\sum_{z \in \mathbb{Z}} z e^{\rho z} P[Z_1 =z] = 0$. Then we
let $\hat{Z}_1,\hat{Z}_2,...$ be i.i.d. random variables with
their distribution defined in the following fashion:
$$P[\hat{Z}_j =z] \equiv \frac{e^{\rho z}}{C_{\rho}} P[Z_j =z]$$
for all $j$ and $z \in \mathbb{R}$, where $C_{\rho} \equiv \sum_{z
\in \mathbb{Z}} P[\hat{Z}_1 =z]= \sum_{z \in \mathbb{Z}} e^{\rho
z} P[Z_1 =z]$. Then the law of $Z_1,...,Z_n$ conditioned on
$Z_1+...+Z_n=0$ is the same as that of $\hat{Z}_1,...,\hat{Z}_n$
conditioned on $\hat{Z}_1+...+\hat{Z}_n=0$, and the case is
reduced to that of $\mu =0$ as $\bold{E} \hat{Z}_j =0$. We also
estimate the covariance equal to $\hat{C} s(1-t)$ for all $0 \leq
s \leq t \leq 1$, where as before
$$\hat{C} = \lim_{n \rightarrow +\infty} \bold{E}[Z_1^2 \mbox{  }|
\mbox{  } Z_1 +...+ Z_n =0].$$
\end{proof}

Observe that the result can be modified for $X_1, X_2,...$ defined
on a multidimensional lattice  $\mathbb{L} \subset
\mathbb{R}^d$,$d>1$, if we condition on $X_n(1)=\bold{a}(n) =
\bold{a+o(1)} \in \{ z\sqrt{n} \mbox{ : } z \in \bigoplus^n_1
\mathbb{L} \}$. We again let point zero be inside the closed
convex hull of $\{ z \mbox{ : } P[Z_1 = z]>0 \}$. In this case the
process $\tilde{X}_n(t) = X_n(t) + (\bold{a} - \bold{a}(n))t$
converges to the Brownian Bridge $B^{0, \bold{a}}$,
and convergence is uniform whenever
 $\bold{a}(n)$  uniformly converges to zero thanks to
the Local CLT.

\begin{thm}
$\tilde{X}_n(t)$ conditioned on $X_n(1)=\bold{a}(n) =
\bold{a}+o(1)$
 converges weakly
 to the Brownian Bridge.
\end{thm}

Here, as before, if we take $t \in \frac{1}{n} \mathbb{Z} \cap
[0,1]$ and let $\alpha =\frac{k}{\sqrt{n}}$, then

\begin{eqnarray*}
P[X_n(t)=\alpha \mbox{  }|\mbox{  } X_n(1)=\bold{a}(n)] & = &
 \frac{(\frac{1}{\sqrt{n}} \Phi_{\sigma \sqrt{t}}(\alpha)
 + o(\frac{1}{\sqrt{n}}))
 (\frac{1}{\sqrt{n}} \Phi_{\sigma \sqrt{1-t}}(\bold{a}(n) -\alpha)
 + o(\frac{1}{\sqrt{n}}))}
 {\frac{1}{\sqrt{n}} \Phi_{\sigma}(\bold{a}(n)) + o(\frac{1}{\sqrt{n}})}
 \\
 & = & \frac{1}{\sqrt{n}} \Phi_{\sigma \sqrt{t(1-t)}}
 (\alpha-\bold{a}(n)t) + o(\frac{1}{\sqrt{n}}) .
\end{eqnarray*}

\subsection{General Case.}\label{bb:gen}

 As before, for a given non-zero vector $\bold{\vec{a}} \in \mathbb{Z}^d$,
 we let $X_1, X_2,...$ be i.i.d. random variables on $\mathbb{Z}^d$ with
 the span of the lattice distribution equal to one (see \cite{durrett})
 such that the probability $P[\bold{\vec{a}} \cdot X_1 >0] =1$,
 the mean $\mu = \bold{E}X_1 <\infty$ and there is a constant
 $\bar{\lambda} >0$ such that the moment-generating function
 $$\bold{E}(e^{\theta \cdot X_1}) <\infty$$
 for all $\theta \in B_{\bar{\lambda}}$. Also we let
 $\bold{P_{\vec{a}}}$ denote the projection
 map on $<\bold{\vec{a}}>$ and
 $\bold{P^{\bot}_{\vec{a}}}$
 denote the orthogonal projection on $<\bold{\vec{a}}>^{\bot}$. Now we can decompose
 the mean $\mu = \mu_{a} \times \mu_{or}$, where
 $\mu_{a} \equiv \bold{P_{\vec{a}}}\mu$ and
 $\mu_{or} \equiv \bold{P^{\bot}_{\vec{a}}}\mu$.\\

 As before we introduce a new basis $\{ \vec{f_1},\vec{f_2},..., \vec{f_d} \}$,
 where $\vec{f_1} = \frac{\bold{\vec{a}}}{\| \bold{\vec{a}} \|}$. We again use
 $[\cdot,\cdot]_f \in \mathbb{R} \times \mathbb{R}^{d-1}$ to denote
 the coordinates of a vector with respect to the new basis.
 We denote  $X_i=[T_i,Z_i]_f \in \mathbb{Z} \times \mathbb{Z}^{d-1}$, where
 $[T_i, 0]_f = \bold{P_{\vec{a}}}X_i$ and
 $[0, Z_i]_f = \bold{P^{\bot}_{\vec{a}}}X_i$, and we let
 $X_1+...+X_i =[t_i, Y_i]_f \in \mathbb{Z} \times \mathbb{Z}^{d-1}$.
 Note: $T_i$ and $Z_i$ don't
have to be independent. Interpolating $Y_i$, we get
$$Y(t)=Y_{[t]}+(t-[t])(Y_{[t]+1} - Y_{[t]})$$
for $0 \leq t \leq \infty$ and if we now define
$Y_n(t)\equiv{\frac{Y(nt)}{\sqrt{n}}}$
for $0\leq{t}\leq{1}$, then the following theorem easily follows from the
previous result:

\begin{cor}
$Y_n(t)$ conditioned on $Y_n(1)=0$ converges weakly to the Brownian Bridge.
\end{cor}

 Since  the first coordinate $T_i$ is positive
with probability one, the next step will be to interpolate
$[t_i,Y_i]_f$, and prove that if scaled and conditioned on
$[t_n,Y_n]_f= X_1+...+X_n= [n \| \bold{\vec{a}} \|,0]_f =
n\bold{\vec{a}}$ it will
 converge weakly to the Brownian
Bridge (with the first coordinate being the time axis). Now, the
last theorem implies the result for $P[[T_i,0]_f=\mu_a]=1$, we
want the same result
for $\bold{E}T_i=\| \mu_a \|$ and $VarT_i<\infty$.
\\
We first let $\bar{X}_i \equiv X_i - \mu_a$, then
$\bold{E}\bar{X}_i= \mu_{or}$ and $Var \bar{X}_i <\infty$. We
again interpolate:
$$\bar{X}(t)=\bar{X}_{[t]}+(t-[t])(\bar{X}_{[t]+1} - \bar{X}_{[t]})$$
for $0 \leq t \leq \infty$, and scale
 $\bar{X}_k(t)\equiv{\frac{\bar{X}(kt)}{\sqrt{k}}}$.
Note: the last $d-1$ coordinates of $\bar{X}_k(t)$ w.r.t. the new
basis are $Y_k(t)$ (e.g.
 $\bold{P^{\bot}_{\vec{a}}}\bar{X}_k(t) = [0, Y_k(t)]_f$).
 \\
  From here on we denote $S_j \equiv [t_j, Y_j]_f =X_1 +...+X_j$
and $\bar{S}_j \equiv \bar{X}_1 +...+\bar{X}_j =S_j- j\mu_a$ for
any positive integer $j$. As a first important step, we state
another important
\begin{cor}
For $k=k(n)= [\frac{n \| \bold{\vec{a}} \|}{\| \mu_a \|} + k_0
\sqrt{n}]$,
 $\{ \bar{X}_k(t) -(k_0\sqrt{\frac{\| \mu_a \|}{\|
\bold{\vec{a}} \|}} \mu_a +\frac{n\bold{\vec{a}} -
k\mu_a}{\sqrt{k}})t \}$ conditioned on
 $\bar{X}_k(1)= n \bold{\vec{a}} -k\mu_a$
 ${(e.g. [t_k,Y_k]_f= n \bold{\vec{a}} )}$
converges weakly to the Brownian Bridge $B^{0, - k_0
\sqrt{\frac{\| \mu_a \|}{\| \bold{\vec{a}} \|}} \mu_a}$.
\end{cor}
Observe that $n \bold{\vec{a}} -k\mu_a  = -k_0 \sqrt{n}\mu_a +o(\sqrt{n})$
and that the convergence is uniform for all $k_0$ in a
compact set . Now, looking only at the last $d-1$ coordinates of
$\bar{X}_k(t)$, w.r.t. the new basis the last Corollary implies:

\begin{lem}
For $k=k(n)= [\frac{n \| \bold{\vec{a}} \|}{\| \mu_a \|} + k_0
\sqrt{n}]$, $Y_k(t)$ conditioned on $t_k=n \| \bold{\vec{a}} \|$
and $Y_k(1)=0$ converges weakly to the Brownian Bridge.
\end{lem}
Note that convergence is uniform for $k_0$ in a compact set.
\\
What the Lemma above says is the following: the interpolation of
$[\frac{i}{k}, \frac{1}{\sqrt{k}}Y_i]_f$ conditioned on $[t_k,
Y_k]_f=n \bold{\vec{a}}$ converges to Time$\times$Brownian Bridge.
Now, define the process $[t, Y_{n,k}^*(t)]_f$ to be the
interpolation of $[\frac{1}{n \| \bold{\vec{a}} \|}t_i,
\frac{1}{\sqrt{n}}Y_i]^{i=0,1,...,k}_f$, then

\begin{thm}
For $k=k(n)= [\frac{n \| \bold{\vec{a}} \|}{\| \mu_a \|} + k_0
\sqrt{n}]$, $\sqrt{\frac{n}{k}} Y^*_{n,k}(t)$ conditioned on
 $t_k=n \| \bold{\vec{a}} \|$ and $Y_k(1)=0$ converges weakly
to the Brownian Bridge.
\end{thm}

\begin{proof}[Proof:]


Here we observe that the mean $\bold{E} [\frac{t_i}{n \|
\bold{\vec{a}} \|} - \frac{t_{i-1}}{n\| \bold{\vec{a}} \|}]$ is
actually equal to $\frac{\| \mu_a \|}{n \| \bold{\vec{a}} \|}
= \frac{1}{k -k_0 \sqrt{n}}+o(\frac{1}{n})$, and that for a given $\epsilon
>0$, the probability of the $\|[\frac{1}{n \| \bold{\vec{a}}
\|}t_i, \frac{1}{\sqrt{n}}Y_i]_f - [\frac{i}{k},
\frac{1}{\sqrt{k}}Y_i]_f \| = |\frac{t_j}{n \| \bold{\vec{a}} \|}
- \frac{j}{k}|$ exceeding $\epsilon$ for some $j\leq k$,

\begin{eqnarray*}
P[ \max_{0\leq j \leq k} |t_j- \frac{n \| \bold{\vec{a}} \|}{k}
j|\geq n\epsilon
  \mbox{  }| \mbox{  } S_n = n \bold{\vec{a}}]
 & \leq &
  P[ \max_{0\leq j \leq k} \| S_j -\frac{n \| \bold{\vec{a}} \| j}{k}\mu_a \|
  \geq  n \epsilon
    \mbox{  } | \mbox{  } S_k = n \bold{\vec{a}}]\\
 & \leq &
  P[ \max_{0\leq j \leq k} |\bar{S}_j| \geq
    n \frac{\epsilon}{2}
   \mbox{  } | \mbox{  } \bar{S}_k =
     [n \| \bold{\vec{a}} \| -k\| \mu_a \|, 0]_f ]\\
 & \rightarrow & 0
\end{eqnarray*}

as $n \rightarrow +\infty$ since
 $n \| \bold{\vec{a}} \| -k\| \mu_a \|
  = -\| \mu_a \| k_0 \sqrt{n} + o(\sqrt{n})$.

\end{proof}

Now, the next step is to prove that the process
$$\{Y^*_{n,k} \mbox{ for some } k \mbox{ such that } [t_k, Y_k]_f= n \bold{\vec{a}} \}$$
conditioned on the existence of such $k$ converges  weakly to the
Brownian Bridge.\\

First of all the last theorem implies
\begin{lem}
For given $k=k(n)= [\frac{n \| \bold{\vec{a}} \|}{\| \mu_a \|} +
k_0 \sqrt{n}]$, $Y^*_{n,k}(t)$ conditioned on $t_k=n \|
\bold{\vec{a}} \|$ and $Y_k(1)=0$ converges weakly to the Brownian
Bridge.
\end{lem}

For a fixed $M>0$, convergence is also uniform on $k \in [\frac{n
\| \bold{\vec{a}} \|}{\| \mu_a \|} -M\sqrt{n},
               \frac{n \| \bold{\vec{a}} \|}{\| \mu_a \|} +M\sqrt{n}]$.
For the future purposes we denote $\kappa \equiv \frac{\| \mu_a
\|}{\| \bold{\vec{a}} \|}$ and
$I_M \equiv [\frac{n}{\kappa} -M\sqrt{n}, \frac{n}{\kappa}+M\sqrt{n}] \bigcap
 \mathbb{Z}$.\\

Finally, we want to prove the following technical result, in which
we use the uniformity of convergence for all $k=k(n) \in I_M$ and
the truncation techniques to show the convergence of $Y^*_{n,k}$
to the Brownian Bridge in case when we condition only on the
existence of such $k$.
\begin{TECHthm}
The process
$$\{ Y^*_{n,k}\mbox{ for some } k \mbox{ such that } [t_k, Y_k]_f = n \bold{\vec{a}} \}$$
conditioned on the existence of such $k$ converges  weakly to the
Brownian Bridge.
\end{TECHthm}

\begin{proof}[Proof:]
Take $M$ large, notice that for $A \subset C^{d-1}[0,1]$,

$$\max_{k \in I_M}
  |P[Y^*_k \in A  \mid [t_k, Y_k]_f = n \bold{\vec{a}}]
  -P[ B^o \in A] | = o(1),$$

where the Brownian Bridge $B^o$ is scaled up to the same constant
for all
those $k$.\\

Hence,
$$lim_{n \rightarrow +\infty}
  \frac{\sum_{k \in I_M}
  P[S_k= n \bold{\vec{a}}]P[Y^*_{n,k} \in A | S_k= n \bold{\vec{a}} ]}
  {\sum_{k \in I_M} P[S_k= n \bold{\vec{a}} ]} = P[ B^o \in A].$$

Therefore we are only left to prove the truncation argument as $M
\rightarrow +\infty$. Now, for any $\epsilon
>0$ there exists $M>0$ such that
$$(1+\epsilon) \sum_{k \in I_M} P[S_k= n \bold{\vec{a}} ]
 \leq \sum_k   P[S_k= n \bold{\vec{a}} ]
 \leq (1+2\epsilon) \sum_{k \in I_M} P[S_k= n \bold{\vec{a}} ]$$
for $n$ large enough, as by the large deviation upper bound, there
is a constant $\bar{C}_{LD} >0$ such that
$$P[S_k= n \bold{\vec{a}} ] \leq
  e^{- \bar{C}_{LD} \frac{(n-k\kappa)^2}{k} \wedge |n-k\kappa|} ,$$
and therefore $\exists C_{LD}>0$ such that
$$\sum_{|n-k \kappa| >n^{2/3} }
  P[S_k= n \bold{\vec{a}}] < e^{-C_{LD} n^{1/3} }.$$
 Also, by the local CLT,
$$P[S_k= n \bold{\vec{a}}]=P[\bar{S}_k= (n-k \kappa) \bold{\vec{a}}]
 =\frac{1}{k^{d/2} \sqrt{Var\bar{X}_1 (2 \pi)^d }}
 e^{-\frac{1}{2Var\bar{X}_1} \frac{(n-k \kappa)^2}{k} }
 +o(\frac{1}{k^{d/2}})$$
implying
 $$\sum_{|n-k \kappa| \leq n^{2/3} }
  P[S_k= n \bold{\vec{a}} ]
 = \frac{1}{ n^{\frac{d-1}{2}} }  [\int_{-\infty}^{+\infty}
 \frac{1}{\sqrt{Var\bar{X}_1 (2 \pi)^d }}
 e^{-\frac{x^2}{2Var\bar{X}_1} } dx +o(1)]$$
where
 $$\sum_{k \in I_M} P[S_k= n \bold{\vec{a}} ]
 = \frac{1}{ n^{\frac{d-1}{2}} }  [\int_{-M}^M
  \frac{1}{\sqrt{Var\bar{X}_1 (2 \pi)^d }}
 e^{-\frac{x^2}{2Var\bar{X}_1} } dx +o(1)].$$
Therefore

$$\frac{1}{1+2\epsilon}  \frac{\sum_{k \in I_M}
  P[S_k= n \bold{\vec{a}}]P[Y^*_{n,k} \in A | S_k= n \bold{\vec{a}}]}
  {\sum_{k \in I_M} P[S_k=n \bold{\vec{a}}]}
  \leq   \frac{\sum_{k}
  P[S_k= n \bold{\vec{a}} ]P[Y^*_{n,k} \in A | S_k= n \bold{\vec{a}} ]}
  {\sum_{k} P[S_k= n \bold{\vec{a}} ]}$$

$$ \leq \frac{1}{1+\epsilon}  \frac{\sum_{k \in I_M}
  P[S_k= n \bold{\vec{a}} ]P[Y^*_{n,k} \in A | S_k= n \bold{\vec{a}} ]}
  {\sum_{k \in I_M} P[S_k= n \bold{\vec{a}} ]}$$

for all $A \subset C^{d-1}[0,1]$. Taking the $\liminf$ and
$\limsup$ of the fraction in the middle completes the proof.

\end{proof}

\section*{Acknowledgements}

 The author wishes to thank D.Ioffe, who posed the problem, and A.Dembo
 for providing him with valuable and insightful comments and suggestions
 concerning the matter of this research.

\bibliographystyle{amsplain}

\end{document}